\documentclass{amsproc}
\usepackage{breqn,float,amssymb,pdflscape,rotating,setspace}
\makeatletter
\@namedef{subjclassname@2020}{%
  \textup{2020} Mathematics Subject Classification}
\makeatother
\newtheorem{theorem}{Theorem}[section]

\theoremstyle{definition}

\newtheorem{example}[theorem]{Example}

\theoremstyle{remark}

\numberwithin{equation}{section}
\allowdisplaybreaks
\begin{document}

\title{Series involving rational, factorial and power functions}


\author{Robert Reynolds}
\address[Robert Reynolds]{Department of Mathematics and Statistics, York University, Toronto, ON, Canada, M3J1P3}
\email[Corresponding author]{milver73@gmail.com}
\thanks{}




\date{}

\dedicatory{}

\begin{abstract}
This is an anthology of series involving rational, factorial, and power functions expressed in terms of special functions. New finite expansions involving quotient functions expressed in terms of the Hurwitz-Lerch zeta function are given. These results represent a new form of expressing this special function as a finite series where contour integration is required for derivation. Extended series previously known and derived are extended using differential equations and algebraic methods.
\end{abstract}
\maketitle
%
%
\section{Introduction}
In both pure research and applied sciences, one often obtains a finite or infinite series and needs to know if there exists an equivalent closed form. A compilation of finite and/or infinite series can serve two essential purposes. First, it can provide the means to express some given mathematical function as a sum of more elementary functions. The second use for a compilation of series is a reverse of this process. Such a list can be used as a means for learning which function is the closed form for a given series.
The book of Hansen \cite{hansen}, contains an exhaustive list of series and products  used widely in physical sciences. In this work we focus on extending Hansen's series involving rational, factorial, and power functions. These series extensions are achieved using contour integration, first order differential equations, and algebraic methods. The derivations of these series are expressed in terms of the Hurwitz-Lerch zeta function. Literature with finite series of quotient functions expressed in terms of the Hurwitz-Lerch zeta function is not widely available to the best of our knowledge. The aim of this current work to assist in bridging this gap by providing a concise table with summarized derivations on this topic. The Hurwitz-Lerch zeta function is used in the study of distribution of prime numbers see (25.16.1) in \cite{dlmf} and Euler sums, see (25.16) in \cite{dlmf} to name a few. Physical applications include analogies between the distribution of the zeros of $\zeta(s)$ on the critical line and of semiclassical quantum eigenvalues. This relates to a suggestion of Hilbert and P\"{o}lya that the zeros are eigenvalues of some operator, and the Riemann hypothesis is true if that operator is Hermitian, see (25.17) in \cite{dlmf}.
\section{Table of double series}
In this section we apply the method in \cite{reyn4}, to the quoted equations to derive the finite double series in terms of the Hurwitz-Lerch zeta function. These forms are useful in simplifying definite integrals written in terms of the Hurwitz-Lerch zeta function and thus simplified in terms of fundamental constants whenever possible. Previous work in this area was conducted by Guillera and Sondow, \cite{guillera}. These forms can also be considered as  recurrence identities involving distant neighbours. These forms can also be considered as elementary translation properties for fixed positive integers $n$, see \cite{bradley}.
\begin{example}
From equation (10.33.16) in \cite{hansen}
\begin{multline}\label{103316}
\sum _{p=0}^n \sum _{j=0}^{2 p} \frac{(-1)^{j+p} e^{i m (j-p)} \binom{2 p}{j} (-n)_p (1+n)_p}{(a+j-p)^k (2
   p)!}\\=e^{-i m n} \Phi \left(-e^{i m},k,a-n\right)+e^{i m (1+n)} \Phi \left(-e^{i m},k,1+a+n\right)
\end{multline}
\end{example}
\begin{example}
From equation (10.33.17) in \cite{hansen}
\begin{multline}
\sum _{p=0}^n \sum _{j=0}^{2 p} \frac{e^{(-j+p) z} (a+j-p)^{-s} \binom{2 p}{j} (-n)_p (1+n)_p}{(2 p)!}\\
=(-1)^n
   e^{n z} \left(\Phi \left(e^{-z},s,a-n\right)-e^{-((1+2 n) z)} \Phi \left(e^{-z},s,1+a+n\right)\right)
\end{multline}
where $a\neq \mathbb{Z}$
\end{example}
\begin{example}
From equation (10.33.18) in \cite{hansen}
\begin{multline}
\sum _{p=0}^n \sum _{j=0}^{2 p} \frac{e^{(j-p) z} (a+j-p)^s \binom{2 p}{j} (-n)_p (1+n)_p}{(2 p)!}\\
=(-1)^{n+s}
   e^{-((1+n) z)} \left(e^{z+2 n z} \Phi \left(e^{-z},-s,-a-n\right)-\Phi \left(e^{-z},-s,1-a+n\right)\right)
\end{multline}
where $a\neq \mathbb{Z},Im(a)\geq 0$
\end{example}
\begin{example}
From equation (10.33.20) in \cite{hansen}
\begin{multline}
\sum _{p=0}^n \sum _{j=0}^{2 p} \frac{(-1)^p e^{(-j+p) z} (a+j-p)^{-s} \binom{2 p}{j} (-n)_p (1+n)_p}{(2
   p)!}\\
=(-1)^{-s} e^{-n z} \left(\Phi \left(-e^z,s,-a-n\right)+e^{z+2 n z} \Phi
   \left(-e^z,s,1-a+n\right)\right)
\end{multline}
where $a\neq \mathbb{Z},Im(a)\geq 0$
\end{example}
\begin{example}
From equation (10.34.15) in \cite{hansen}
\begin{multline}
\sum _{p=0}^{n-1} \sum _{j=0}^{2 p+1} \frac{(-1)^{j+p} e^{i (j-p) z} (a+j-p)^{-s} \binom{1+2 p}{j} (1-n)_p (1+n)_p}{\Gamma (2+2
   p)}\\
\frac{e^{-i (-1+n) z} \left(\Phi \left(-e^{i z},s,1+a-n\right)-e^{2 i n z} \Phi \left(-e^{i
   z},s,1+a+n\right)\right)}{n}
\end{multline}
\end{example}
\begin{example}
From equation (10.34.16) in \cite{hansen}
\begin{multline}
\sum _{p=0}^{n-1} \sum _{j=0}^{2 p+1} \frac{e^{(-j+p) z} (a+j-p)^{-s} \binom{1+2 p}{j} (1-n)_p (1+n)_p}{\Gamma (2+2
   p)}\\
=\frac{(-1)^{1+n} e^{-((1+n) z)} \left(e^{2 n z} \Phi \left(e^{-z},s,1+a-n\right)-\Phi
   \left(e^{-z},s,1+a+n\right)\right)}{n}
\end{multline}
where $Re(a)\leq 0$
\end{example}
\begin{example}
From equation (10.34.17) in \cite{hansen}
\begin{multline}
\sum _{p=0}^{n-1} \sum _{j=0}^{2 p+1} \frac{e^{(j-p) z} (a+j-p)^{-s} \binom{1+2 p}{j} (1-n)_p (1+n)_p}{\Gamma (2+2 p)}\\
=\frac{(-1)^{1+n-s} e^{-n z} \left(e^{2 n z} \Phi \left(e^{-z},s,-a-n\right)-\Phi
   \left(e^{-z},s,-a+n\right)\right)}{n}
\end{multline}
where $Re(a)\leq 2\pi$
\end{example}
\begin{example}
From equation (10.34.20) in \cite{hansen}
\begin{multline}
\sum _{p=0}^{n-1} \sum _{j=0}^{2 p+1} \frac{(-1)^{j+p} e^{(-j+p) z} (a-j+p)^{-s} \binom{1+2 p}{j} (1-n)_p (1+n)_p}{\Gamma (2+2 p)}\\
=\frac{(-1)^{-s} e^{-((1+n) z)} \left(e^{2 n z} \Phi
   \left(-e^{-z},s,1-a-n\right)-\Phi \left(-e^{-z},s,1-a+n\right)\right)}{n}
\end{multline}
where $Re(a)\leq -\pi$
\end{example}
\begin{example}
From equation (10.27.7) in \cite{hansen}
\begin{multline}
\sum _{p=0}^n \sum _{j=0}^{2 (n-p)} \frac{e^{i (j+3 p) \pi -(j+p) z} (a+j-n+p)^{-s} \binom{2 n-2 p}{j} (-2 n)_{2 p}}{p! (-2 n)_p}\\
=e^{2 i n \pi } \left(\Phi \left(-e^{-z},s,a-n\right)+e^{-((1+2 n) z)} \Phi
   \left(-e^{-z},s,1+a+n\right)\right)
\end{multline}
\end{example}
\begin{example}
From equation (10.27.8) in \cite{hansen}
\begin{multline}\label{10278}
\sum _{p=0}^{n-1} \sum _{j=0}^{2 n-2 p-1} \frac{e^{i (j+3 p) \pi +(j+p) z} \binom{2 n-2 p-1}{j} (1-2 n)_{2
   p}}{p! (1-2 n)_p (a+j-n+p)^s}\\
=\Phi \left(-e^z,s,a-n\right)-e^{2 n z} \Phi \left(-e^z,s,a+n\right)
\end{multline}
where $a\neq \mathbb{Z}$
\end{example}
\begin{example}
From equation (10.27.9) in \cite{hansen}
\begin{multline}
\sum _{p=0}^{\left\lfloor \frac{n}{2}\right\rfloor } \sum _{j=0}^{n-2 p} \frac{e^{-((j+p) z)} (a+j+p)^{-s} \binom{n-2 p}{j} (-n)_{2 p}}{p! (-n)_p}\\
=\Phi \left(e^{-z},s,a\right)-e^{-((1+n) z)} \Phi
   \left(e^{-z},s,1+a+n\right)
\end{multline}
\end{example}
\begin{example}
From equation (10.27.10) in \cite{hansen}
\begin{multline}
\sum _{p=0}^{\left\lfloor \frac{n}{2}\right\rfloor } \sum _{j=0}^{n-2 p} \frac{e^{(j+p) z} (a+j+p)^{-s} \binom{n-2 p}{j} (-n)_{2 p}}{p! (-n)_p}\\
=(-1)^{-s} e^{-z} \left(-\Phi
   \left(e^{-z},s,1-a\right)+e^{(1+n) z} \Phi \left(e^{-z},s,-a-n\right)\right)
\end{multline}
where $Re(a)\leq -\pi$
\end{example}
\begin{example}
From equation (10.27.12) in \cite{hansen}
\begin{multline}
\sum _{p=0}^n \sum _{j=0}^{2 n-2 p} \frac{(-1)^{j+p} e^{-((j+p) z)} (a-j-p)^{-s} \binom{2 n-2 p}{j} (-2 n)_{2 p}}{p! (-2 n)_p}\\
=(-1)^{-s} \left(\Phi \left(-e^{-z},s,-a\right)+e^{-((1+2 n) z)} \Phi
   \left(-e^{-z},s,1-a+2 n\right)\right)
\end{multline}
where $Re(a)\leq -\pi$
\end{example}
%
%
%
\section{Extended series involving rational, factorial and power functions.}
In this section we apply the method in section (7) in \cite{reyn_plos} to the stated equations to derive the corresponding examples.
\begin{example}
From equation (10.49.26) in \cite{hansen}. Extended McClintoch series.
\begin{multline}
\sum _{k=1}^{\infty } \frac{\left((-1+y) y^{-1-a}\right)^k \Gamma ((1+a) k)}{(k+1)! \Gamma (1+a k)}\\
=\frac{-a+y+a y-y^{1+a}-a \log (y)-a^2 \log (y)+a y \log (y)+a^2 y \log (y)}{a (1+a) (-1+y)}
\end{multline}
where $|Re(a)| < Re(y)$
\end{example}
\begin{example}
\begin{multline}
\sum _{k=1}^{\infty } \frac{\left((-1+y) y^{-1-a}\right)^k \Gamma ((1+a) k)}{(k+2)! \Gamma (1+a k)}
=\frac{y^{2 (1+a)}}{4 a (1+a) (1+2 a) (1-y)^2}\\+\frac{1}{4 a (1+a) (1+2 a) (1-y)^2}\left(6 a^2 (-1+y)^2-y \left(4-3 y-4 y^a+4 y^{1+a}\right)\right. \\ \left.
-a
   (-1+y) \left(3-9 y+8 y^{1+a}\right)+2 a \left(1+3 a+2 a^2\right) (-1+y)^2 \log (y)\right)
\end{multline}
where $|Re(a)| < Re(y)$\end{example}
\begin{example}
From equation (10.49.27) in \cite{hansen}. Extended P\'{o}lya-Szeg\"{o} series.
\begin{multline}
\sum _{k=0}^{\infty } \frac{\left((-1+y) y^{-1-a}\right)^k (-1+b)_{(1+a) k}}{(k+1)! (b)_{a k}}\\
=-\frac{(1-b) y^{1+a}}{(1+a-b) (2+a-b) (1-y)}+\frac{y^{-1+b} \left(\frac{1+a}{2+a-b}-\frac{a
   y}{1+a-b}\right)}{1-y}
\end{multline}
where $|Re(a)| < Re(y)$
\end{example}
\begin{example}
From equation (10.49.28) in \cite{hansen}
\begin{multline}
\sum _{k=0}^{\infty } \frac{\left((-1+y) y^{-1-a}\right)^k (b)_{(1+a) k}}{(k+1)! (b)_{a k}}=-\frac{y^{1+a}}{(1+a-b) (1-y)}-\frac{y^b}{(-1-a+b) (1-y)}
\end{multline}
where $|Re(a)| < Re(y)$
\end{example}
\begin{example}
From equation (10.49.29) in \cite{hansen}
\begin{multline}
\sum _{k=0}^{\infty } \frac{\left((-1+y) y^{-1-a}\right)^k (1+b)_{(1+a) k}}{(k+1)! (b)_{a k}}=\frac{y^{1+a}}{b (1-y)}+\frac{y^{1+b}}{b (-1+a (-1+y)) (1-y)}
\end{multline}
where $|Re(a)| < Re(y)$
\end{example}
\section{Definite integrals in terms special cases of the Hurwitz-Lerch zeta function}
In this section we look at examples using the above Hurwitz-Lerch zeta function double series to simplify definite integrals.
\begin{example}
From equation (7) in \cite{reyn_log} and equation (\ref{10278}).
\begin{multline}
\int_0^{\infty } \frac{x^{3/2} \log (-x) \log (\log (-x))}{1+x^5} \, dx=-\frac{1}{25} \pi ^2 \left(\frac{5 \pi
   }{2}+i \log \left(\frac{3125 A^{12}}{16384 \sqrt[3]{2} e \pi ^5}\right)\right)\
\end{multline}
where $A$ is Glaisher's constant given in equation (5.17.6) in \cite{dlmf}.
\end{example}
\begin{example}
From equation (7) in \cite{reyn_log} and equation (\ref{10278}).
\begin{equation}
\int_0^{\infty } \frac{x^4 \log ^2(i x) \log (\log (i x))}{1+x^{10}} \, dx=\frac{3}{250} \pi ^3 \left(i \pi
   +\log \left(\frac{25 e^{\frac{7 \zeta (3)}{6 \pi ^2}}}{4\times 2^{2/3} \pi ^2}\right)\right)
\end{equation}
where $\zeta(3)$ is Ap\'{e}ry's constant.
\end{example}
\begin{example}
From equation (7) in \cite{reyn_log} and equation (\ref{10278}).
\begin{equation}
\int_0^{\infty } \frac{x \log (-x) \log (\log (-x))}{1+x^4} \, dx=-\frac{1}{8} \pi  \left(2 i C+\pi  \left(\pi
   -i \log \left(\frac{27 \pi ^2}{16}\right)\right)\right)
\end{equation}
where $C$ is Catalan's constant.
\end{example}
\begin{example}
From equation (\ref{103316})
\begin{multline}
\sum _{p=0}^n \sum _{j=0}^{2 p} \frac{(-1)^{j+p} e^{i m (j-p)}
   (a+j-p)^{-k} \binom{2 p}{j} \,
   _2F_1\left(k,\frac{1}{u};1+\frac{1}{u};-\frac{1}{a+j-p}\right) (-n)_p
   (1+n)_p}{(2 p)!}\\
=\int_0^1 \left(e^{-i m n} \Phi \left(-e^{i
   m},k,a-n+x^u\right)+e^{i m (1+n)} \Phi \left(-e^{i
   m},k,1+a+n+x^u\right)\right) \, dx
\end{multline}
where $Re(k)>0,Re(u)>Re(n),Re(a)>\pi$.
\end{example}
%
%
%

%
\end{document}